\documentclass{amsart}
\usepackage{amssymb,color,cancel,mathrsfs, manfnt
}

\def\sF{\mathscr F}

\def\fa{\mathfrak{a}}

\def\fn{\mathfrak{n}} 
\def\Z{{\mathbb Z}}   \def\N{{\mathrm N}} \def\Q{{\mathbb Q}}

\def \ring {\Z[\omega]}

\DeclareMathOperator{\tr}{Tr}

\DeclareMathOperator{\re}{Re}

\def \e{e}
\def\cond#1{n_{#1}}

\renewcommand{\mod}[1]{\,{\rm mod}\,#1}

\def\su#1{\sum_{\substack{#1}}}
\def\ps#1{\sideset{}{^*}\sum_{\substack{#1}}}
\def\pr#1{\prod_{\substack{#1}}}

\def\bs#1{\begin{equation*} \begin{split} #1 \end{split} \end{equation*}}
\def\bsc#1{\begin{equation} \begin{split} #1 \end{split} \end{equation}}
\def\eqs#1{\begin{equation*} #1 \end{equation*}}
\def\eqn#1{\begin{equation} #1 \end{equation}}
\def\mult#1{\begin{multline*}#1\end{multline*}}
\def\multn#1{\begin{multline}#1\end{multline}}

\def\({\left(} \def\){\right)} \def\[{\left[} \def\]{\right]} 
 
\def\le{\leqslant} \def\ge{\geqslant}
\def\eps{{\varepsilon}}

\definecolor{pink}{rgb}{1,.2,.6}
\definecolor{orange}{rgb}{0.7,0.3,0}
\definecolor{blue}{rgb}{.2,.6,.75}
\definecolor{green}{rgb}{.4,.7,.4}
\definecolor{mycolor}{RGB}{210,15,35}
\def\red#1{\textcolor{mycolor}{#1}}

\newcommand*\cube{\mbox{\mancube}}

\newtheorem{lem}{Lemma}[section]

\newtheorem{cor}[lem]{Corollary}

\newtheorem{theorem}[lem]{Theorem}
\theoremstyle{definition}

\newtheorem{rem}[lem]{Remark}

\begin{document}
\author{Ahmet M. G\"ulo\u{g}lu  } 
\address{Department of Mathematics, Bilkent University, Ankara, Turkey} \email{guloglua@fen.bilkent.edu.tr} 


\title[Non-vanishing of Cubic Dirichlet L-Functions]{Non-vanishing of Cubic Dirichlet L-Functions over the Eisenstein Field}

\begin{abstract} 
We establish an asymptotic formula for the first moment and derive an upper bound for the second moment of L-functions associated with the complete family of primitive cubic Dirichlet characters defined over the Eisenstein field. Our results are unconditional, and indicate that there are infinitely many characters within this family for which the L-function $L(s,\chi)$ does not vanish at the central point $s=1/2$.
\end{abstract} 

\maketitle
\section{Introduction}
According to Chowla's conjecture \cite{Chow}, it is stated that for any real non-principal Dirichlet character $\chi$, the value of $L(1/2, \chi)$ cannot be equal to zero. However, Li \cite{Li} demonstrated that there are infinitely many quadratic Dirichlet $L$-functions over the rational function field for which $L(1/2, \chi) = 0$. Nonetheless, it is widely believed that among all quadratic characters, the number of such characters with $L(1/2, \chi) = 0$ should have zero density.

In the case of quadratic Dirichlet $L$-functions, Özlük and Snyder \cite{OS-1999}, assuming the Generalized Riemann Hypothesis (GRH), computed the one-level density for the low-lying zeroes in the family and demonstrated that at least 15/16 of these functions do not vanish at $s=1/2$. The conjectures put forth by Katz and Sarnak \cite{KS-book} suggest that $L(1/2, \chi) \neq 0$ for almost all quadratic Dirichlet $L$-functions. Soundararajan \cite{Sound}, by calculating the first two mollified moments, proved that at least 87.5\% of quadratic Dirichlet $L$-functions do not vanish at $s = 1/2$ without assuming GRH. It is worth noting that using only the first two (non-mollified) moments does not yield a positive proportion of non-vanishing, as their growth rate is too rapid (refer to Conjecture 1.5.3 in \cite{CFKRS} and the work of Jutila \cite{Jut}).

In the function field case, Bui and Florea \cite{BuiFlo2018} computed the one-level density and obtained a proportion of at least 94\% for non-vanishing quadratic Dirichlet $L$-functions.

In this paper, we consider the family of $L$-functions attached to primitive cubic Dirichlet characters defined over the Eisenstein field $K = \Q(\omega)$, where $\omega = e^{2\pi i/3}$. Various researchers have investigated the one-level density for families of cubic characters over both $\Q$ and $K$, employing test functions whose Fourier transforms are limited to the interval $(-1, 1)$. These families possess unitary symmetry. However, to achieve a positive proportion of non-vanishing, it becomes necessary to extend beyond the interval $(-1, 1)$. David and the author accomplished this for a subset of cubic characters over $K$ in their work \cite{DG}. They employed test functions with support in $(-13/11, 13/11)$, which yielded a proportion of non-vanishing equal to $2/13$. Notably, these results are obtained assuming the Generalized Riemann Hypothesis (GRH).

Regarding the moments of $L(1/2, \chi)$, where $\chi$ is a primitive cubic character, Baier and Young \cite{BaYo} computed the first moment over $\Q$, Luo \cite{Luo} calculated it for the same thin sub-family as in \cite{DG} over $K$, and David, Florea, and Lalin \cite{DFL-1} did so over function fields. The first two papers utilized smoothed sums over the respective families. In all three cases, the authors obtained lower bounds for the number of non-vanishing cubic twists but not positive proportions, relying on upper bounds for higher moments. The asymptotic behavior of the second moment for cubic Dirichlet $L$-functions remains an open question for function fields and number fields.

David, Florea, and Lalin \cite{DFL-2} established, using mollified moments of $L$-functions, that there is a positive proportion of non-vanishing among cubic Dirichlet $L$-functions at $s = 1/2$ over function fields in the non-Kummer case. Using the same methodology, Yesilyurt and the author have recently demonstrated a positive proportion of non-vanishing for the entire family of cubic Dirichlet characters over $K$. This result had not previously been established for the complete family of cubic characters. It is important to note, however, that it relies on the GRH for the higher moments. 

In this study, we unconditionally obtain an asymptotic formula for the first moment and an upper bound for the second moment of the $L$-functions associated with the entire family of primitive cubic Dirichlet characters defined over the Eisenstein field $K$. Our approach is similar to that of Luo \cite{Luo}, although we consider characters with conductors of bounded norm rather than employing a smoothed sum over the family.

The result presented in this paper builds upon the methods developed in \cite{DG}, which rely on the profound contributions of Kubota and Patterson regarding the average of cubic Gauss sums. Furthermore, we utilize Heath-Brown's cubic large sieve inequality as described in \cite[Theorem 2]{HB2000}. These techniques form the foundation for our analysis and enable us to derive the following result.

\begin{theorem} \label{thm:moments}
Consider the family $\mathcal{F}(X)$, defined in equation \eqref{family}, consisting of all primitive cubic Dirichlet characters whose conductors have a norm not exceeding $X$. Then,
\eqs{
\sum_{\chi \in \sF(X)} L(1/2, \chi) = D X\log X + EX + O (X^{65/66 + \eps}),}
where $D$ and $E$ are given in \eqref{AsymptoticConstant}. Furthermore, 
\eqn{\label{secondmoment}
\sum_{\chi \in \sF(X)} |L(1/2, \chi)|^2 \ll X^{7/6+\eps} .} 
\end{theorem}
Using theorem \ref{thm:moments} together with Cauchy-Schwarz inequality we obtain the following result.
\begin{cor}
There exist infinitely many primitive cubic Dirichlet characters $\chi$ such that $L(1/2, \chi) \neq 0$. 
More precisely, the number of such characters whose conductor has a norm not exceeding $X$ is $\gg X^{5/6-\eps}$.
\end{cor}

\section{Preliminary}

\subsection{Notation}
$\N(n) \sim N$ means $N < \N(n) \le 2N$. We write $\N\fa$ for the norm of the ideal $\fa$, and $\N(a)$ for the norm of the ideal $a\ring$. We use 
$e(z) = \exp ( 2\pi i z)$ and put $\tr(z) = z + \bar z$. 

The following lemma is used to find optimal bounds in the proofs of several results throughout the paper. 
\begin{lem}[{\cite[Lemma 2.4]{GraKol}}] \label{balancing} 
	Suppose that $$L(H) = \sum_{i=1}^m A_i H^{a_i} + \sum_{j=1}^n B_j H^{-b_j},$$ where
	$A_i, B_j, a_i, b_j$ are positive, and that $H_1 \le H_2$. Then, there is some $H$ with $H_1 \le H \le H_2$ such that
	$$
	L(H) \ll \sum_{i=1}^m \sum_{j=1}^n \left(A_i^{b_j} B_j^{a_i}\right)^{1/(a_i+b_j)} + \sum_{i=1}^m A_i H_1^{a_i} + \sum_{j=1}^n B_j H_2^{-b_j},
	$$
	where the implied constant depends only on $m$ and $n$.
\end{lem}

\subsection{Cubic Characters}
The ring of integers $\ring$ of $K$ has class number one and six units $\left\{ \pm 1, \pm \omega, \pm \omega^2 \right\}$. Each non-trivial principal ideal $\fn$ co-prime to $3$ has a unique generator $n \equiv 1 \mod 3$.

The cubic Dirichlet characters on $\Z[\omega]$ are given by the cubic residue symbols. For each prime $\pi \in \ring$ with $\pi \equiv 1 \mod 3$, there are two primitive characters of conductor $(\pi)$; the cubic residue character $\chi_\pi$ satisfying 
\eqs{\chi_\pi(\alpha) =  \( \frac{\alpha}{\pi}\)_3 \equiv \alpha^{(\N(\pi)-1)/3} \mod \pi\ring,}
and its conjugate $\overline{\chi}_\pi= \chi_\pi^2$. 

In general, for  $n\in\ring$ with $n \equiv 1 \mod 3$, the cubic residue symbol $\chi_n$ is defined multiplicatively using the characters of prime conductor by
\eqs{\chi_n (\alpha) = \left( \frac{\alpha}{n} \right)_3 = \prod_{\pi^{v_\pi} \| n} \chi_\pi (\alpha)^{v_\pi} .}

Such a character $\chi_n$ is primitive when it is a product of characters of distinct prime conductors, i.e. either $\chi_\pi$ or $\overline \chi_\pi = \chi_{\pi}^2 = \chi_{\pi^2}$.
Moreover, $\chi_n$ is a (cubic) \emph{Hecke} character of conductor $n \ring$ if $\chi_n(\omega) = 1$. 
Since
\eqs{\( \frac{\omega}{n}\)_3 = \prod_{\pi \mid n} \omega^{v_\pi(n) (\N(\pi)-1)/3} = \omega^{\sum_{\pi \mid n} v_\pi(n) (\N(\pi)-1)/3} = \omega^{(\N(n) - 1)/3},}
we conclude that a given Dirichlet character $\chi$ is a primitive cubic Hecke character provided that $\chi = \chi_n$, where 
\begin{enumerate}
	\item $n=n_1 n_2^2$, where $n_1, n_2$ are square-free and co-prime, and
	\item  $\N(n) \equiv 1 \mod 9$, or equivalently, $\N(n_1) \equiv \N(n_2) \mod 9$.
\end{enumerate}
In this case, $\chi$ has conductor $n_1 n_2\ring$. 

We recall the cubic reciprocity theorem (cf. \cite[page 114, Theorem 1]{IR}) for cubic characters. Let $m,n \in \ring, m,n \equiv \pm 1 \mod 3.$ Then,
	$$
	\( \frac{m}{n}\)_3 = \( \frac{n}{m}\)_3.
	$$

\subsection{The family $\sF$ of cubic Dirichlet characters}
We shall consider the family $\sF$ of cubic characters $\chi_{c_1}\overline{\chi_{c_2}}$ given by cubic residue symbols, where $c_1, c_2 \equiv 1 \mod 3\in \Z[\omega]$ are square-free and co-prime, and $c_1c_2^2 \equiv 1 \mod 9$. We naturally exclude the case $c_1=c_2=1$. By a slight abuse of notation (dropping the letter $\chi$), we shall write
\bsc{ \label{family}
	\sF &= \Big\{ c_1c_2^2 \in \ring\setminus\{1\} :  \begin{array}{l} 	c_1, c_2 \equiv 1 \mod 3 \text{ both square-free},\\
		(c_1,c_2)=1, \;c_1c_2^2 \equiv 1 \mod 9
	\end{array} \Big\} \\
	&= \Big\{ qd \in \ring\setminus\{1\} :  \begin{array}{l}
		q, d \equiv 1 \mod 3, q \text{ square-free},\\
		d \mid q, \;qd \equiv 1 \mod 9
	\end{array} \Big\}. 
}
Note that for $c=qd \in \sF$, the conductor of the character $\chi_{qd}$ is $q\ring$. We shall write $\cond{c}$ for the norm of the conductor of $\chi_c$, and $c \in \sF(X)$ to mean that $c\in \sF$ and $\cond c \le X$.

\subsection{Cubic Gauss Sums}
For any $n\equiv 1 \mod 3$, the shifted cubic Gauss sum is defined by 
\eqn{ \label{ShiftedGaussSum}
	g(r,n) = \sum_{\alpha \mod n} \chi_n (\alpha) \e \bigl(\tr( r \alpha/n)\bigr).}

\section{Second Moment}
By the approximate functional equation of $L(s,\chi)$ for $s=1/2$ we have (cf., for example, \cite[Theorem 5.3]{IwaKow}) that
\bsc{\label{approxfnceqn}
L(1/2,\chi_c) &= \sum_{r \ge 0} 3^{-r/2} \sum_{a \equiv 1 \mod 3}  \frac{\chi_c (a)}{\N(a)^{1/2}}  V \biggl( \frac {3^{r}\N(a)} Y \biggr)  \\
&\quad +  \frac{W(\chi_c)}{\sqrt{\cond{c}}}
\sum_{r \ge 0} 3^{-r/2} \sum_{a \equiv 1 \mod 3} \frac{\overline{\chi_c (a)}}{\N(a)^{1/2}} V \biggl(  \frac {3^{r}\N(a)Y} {3\cond{c}}  \biggl),}
where
\eqs{
	V (y)  = \frac 1 {2\pi i} \int_{2-i\infty}^{2+i\infty} (2 \pi y)^{-u} \frac{\Gamma(1/2+u)}{\Gamma(1/2)} \frac{du}{u}.
}
\begin{rem} 
For $0< \alpha \le 1/2-\eps$ with $\eps < 1/2$, and $A > 0$, it follows from Stirling's formula for the Gamma function (see, for example, \cite[Eqn. 5.112]{IwaKow}) by shifting the contour to $\re u = -\alpha$ and to $\re u = A$, respectively, that 
\bsc{\label{Vbound}
		y^a \frac{d^aV}{dy^a} (y) &= \delta_{0,a} + O_{\eps, \alpha} ( y^\alpha)  \\
		y^a \frac{d^aV}{dy^a} (y) &\ll_A y^{-A}
}
for any integer $a \ge 0$, where $\delta_{0,a}$ is the Kronecker delta function. 
\end{rem}
We choose $Y = X^{1/2}$ in this section. By Cauchy-Schwarz inequality
\eqs{|L(1/2,\chi_c)|^2  \le \frac {2\sqrt 3} {\sqrt 3 - 1}
\sum_{r \ge 0} 3^{-r/2} \Bigl( \big| \Sigma_1  \big |^2 + \big | \Sigma_2 \big|^2\Bigr)
}
where 
\bs{\Sigma_1 &= \sum_{a \equiv 1 \mod 3}  \frac{\chi_c (a)}{\N(a)^{1/2}}  V \biggl( \frac {3^{r}\N(a)} Y \biggr) \\
\Sigma_2 &= \sum_{a \equiv 1 \mod 3} \frac{\overline{\chi_c (a)}}{\N(a)^{1/2}} V \biggl(  \frac {3^{r}\N(a)Y} {3\cond{c}}  \biggl).
}
Using partial integration and Cauchy-Schwarz inequality again we have that
\bsc{\label{2ndmomentsum1}
\big| \Sigma_1  \big |^2
&\le 
\int_{1^-}^\infty \bigg|\frac {3^r} Y V' \biggl( \frac {3^r z} Y \biggr) \bigg| dz \cdot \int_{1^-}^\infty \bigg|\frac {3^r} Y V' \biggl( \frac {3^r z} Y \biggr) \bigg| \bigg| \su{a \equiv 1 \mod 3\\ \N(a) \le z}  \frac{\chi_c (a)}{\N(a)^{1/2}} \bigg|^2 dz \\
&\ll \log Y \biggl(\int_1^Y z^{-1} + \int_Y^\infty \frac{Y^2}{3^{2r}z^3} \biggr)
\bigg| \su{a \equiv 1 \mod 3\\ \N(a) \le z}  \frac{\chi_c (a)}{\N(a)^{1/2}} \bigg|^2 dz,
}
where we used the bounds $yV'(y) \ll 1$ when $z \le Y= \sqrt X$, and $yV'(y) \ll y^{-2}$ for $z > \sqrt X$.

Writing $a= a_1a_2^2$ with $a_i \equiv 1 \mod 3$ and $a_1$ square-free, and using Cauchy's inequality shows that
\bs{\bigg| \su{a \equiv 1 \mod 3\\ \N(a) \le z}  \frac{\chi_c (a)}{\N(a)^{1/2}} \bigg|^2 
&=
\bigg|\su{a_2 \equiv 1 \mod 3\\ \N(a_2) \le \sqrt z}  \frac{\chi_c (a_2^2)}{\N(a_2)}
\ps{a_1 \equiv 1 \mod 3\\ \N(a_1a_2^2) \le z}  \frac{\chi_c (a_1)}{\N(a_1)^{1/2}}\bigg| \\
& \ll \log z \su{a_2 \equiv 1 \mod 3\\ \N(a_2) \le \sqrt z}  \frac 1 {\N(a_2)} \bigg| \quad \ps{a_1 \equiv 1 \mod 3\\ \N(a_1a_2^2) \le z}  \frac{\chi_c (a_1)}{\N(a_1)^{1/2}}\bigg|^2.
}
Here and in what follows, $\Sigma^\ast$ indicates summation over square-free integers.
Note that if $A(c) \ge 0$ for all $c\in \sF(X)$, removing the conditions $(c_1, c_2)=1$ and $c_1c_2^2 \equiv 1 \mod 9$ on each $c=c_1c_2^2$, it follows that
\bs{
\su{c \in \sF(X)} A(c) 
&\le  \ps{c_1 \equiv 1 \mod 3\\ \N(c_1) \le \sqrt X}\;  \ps{c_2 \equiv 1 \mod 3\\ \N(c_1c_2) \le X} A (c_1c_2^2)  + \ps{c_2 \equiv 1 \mod 3\\ \N(c_2) \le \sqrt X}\; \ps{c_1 \equiv 1 \mod 3\\ \N(c_1c_2) \le X} A (c_1c_2^2).
}
Applying this idea with 
\eqs{A(c) = \bigg| \quad \ps{a_1 \equiv 1 \mod 3\\ \N(a_1a_2^2) \le z}  \frac{\chi_c (a_1)}{\N(a_1)^{1/2}}\bigg|^2}
and using cubic large sieve inequality \cite[Theorem 2]{HB2000}, we see that
\bsc {\label{2ndmomentsieve}
\su{c \in \sF(X)} \bigg| \su{a \equiv 1 \mod 3\\ \N(a) \le z}  \frac{\chi_c (a)}{\N(a)^{1/2}} \bigg|^2 
&\ll X^{1+\eps}z^\eps +  z^{1+\eps} X^{1/2+\eps} + X^{5/6+\eps} z^{2/3+\eps}.
}
Thus, summing \eqref{2ndmomentsum1} over $c$ and using the last estimate \eqref{2ndmomentsieve} gives that
\bs{\sum_c  \big| \Sigma_1  \big |^2  
&\ll X^{1+\eps} Y^{2\eps} + Y^{1+2\eps} X^{1/2+\eps} + X^{5/6+\eps} Y^{2/3+2\eps} \ll X^{7/6+2\eps}.
}
Similarly we have
\bs {
\big| \Sigma_2  \big |^2 
&\ll \log X \biggl(\int_1^{\sqrt X} \frac 1 z + \int_{\sqrt X}^\infty \frac {X^2}{z^3Y^2} \biggr) \bigg| \su{a \equiv 1 \mod 3\\ \N(a) \le z}  \frac{\chi_c (a^2)}{\N(a)^{1/2}} \bigg|^2 dz 
. }
By summing over $c$ and applying \eqref{2ndmomentsieve} once more, we arrive at the same estimate of $X^{7/6+2\epsilon}$.
Therefore, the claim \eqref{secondmoment} in theorem \ref{thm:moments} follows. 

\section{First Moment}
We write the first moment using \eqref{approxfnceqn} as 
\eqs{
\sum_{c \in \sF(X)} L(1/2, \chi_c) = S_1 + S_2 + S_3,}
where
\bs {
	S_1 &= \sum_{r \ge 0} 3^{-r/2} \su{a \equiv 1 \mod 3\\ a = \cube} \frac 1 {\N(a)^{1/2}} V \biggl( \frac {3^{r}\N(a)} Y \biggl) \su{c \in \sF(X)\\ (c,a)=1} 1    \\
	S_2 &= \sum_{r \ge 0} 3^{-r/2} \su{a \equiv 1 \mod 3\\ a \neq \cube} \frac 1 {\N(a)^{1/2}} V \biggl( \frac {3^{r}\N(a)} Y \biggl) \su{c \in \sF(X)\\ (c,a)=1} \chi_{a} (c)  \\
	S_3 &= \sum_{r \ge 0} 3^{-r/2} \sum_{a \equiv 1 \mod 3} \frac 1 {\N(a)^{1/2}} 
	\sum_{c \in \sF(X)}  \frac{\chi_c (a^2) W(\chi_c)}{\sqrt{\cond{c}}} V \biggl( \frac {3^{r}Y\N(a)} {3\cond{c}} \biggr) .   }

\subsection{The main term $S_1$}
First recalling that $V(y) = 1 + O(y^{1/6-\eps})$ (see \eqref{Vbound}) and trivially estimating the sum over $c$ by $X\log X$ we see that
\eqs{
S_1 = \sum_{r \ge 0} 3^{-r/2} \su{a \equiv 1 \mod 3\\ a = \cube} \frac 1 {\N(a)^{1/2}} \su{c \in \sF(X)\\ (c,a)=1} 1  + O (X^{1+\eps} Y^{-1/6}).
}
Then, one can proceed as in \cite[Lemma 3.1]{DG} to get the main terms with an admissible error. But to show that one can get a better error term, we shall continue. 

Next, removing the conditions $(c_1, c_2) =1,c_1c_2^2 \equiv 1 \mod 9$ and $c_1, c_2$ be square-free, we can write
\bs{\su{c \in \sF(X)\\ (c,a)=1} 1  
&=  \frac 1 {h_{(9)}} \sum_{\psi \mod 9} 
\su{d \equiv 1 \mod 3\\(d, a)=1\\ \N(d)^2 \le X} \mu_K(d)  \psi (d^3)
\su{e_1 \equiv 1\mod 3\\ (e_1,ad)=1\\ \N(de_1)^2 \le X} \mu_K (e_1) \psi  (e_1^2)  \\
& \quad \cdot
\su{e_2 \equiv 1\mod 3\\ (e_2, ad)=1\\ \N(de_1e_2)^2 \le X} \mu_K (e_2) \psi  (e_2^4) \su{c_1 \equiv 1 \mod 3\\  (c_1, ad)=1} \psi  (c_1)
\su{c_2 \equiv 1 \mod 3\\ (c_2, ad)=1\\ \N(c_1c_2) \le W}  \psi  (c_2^2),
}
where $W= X/\N(de_1e_2)^2$. Applying Perron's formula, the sums over $c_1, c_2$ then becomes 
\eqs{\int_{\varrho-iT}^{\varrho+iT} W^s G_\psi (s)  L(s,\psi) L(s,\psi^2) \frac{ds} s + O ( W^\eps + W^{1+\eps}T^{-1} ),
}
where $\varrho = 1 + \log (2W), T \in [1,W]$ and 
\eqs{G_\psi(s) = \prod_{\pi \mid ad} \Bigl(1-\frac{\psi(\pi)}{\N(\pi)^s} \Bigr)\Bigl(1-\frac{\psi^2(\pi)}{\N(\pi)^s} \Bigr) 
. 
}
Moving the contour to $\sigma = \eps$ we pick up the residue at $s=1$ coming from the double pole of $L(s,\psi)L(s,\psi^2)$ when $\psi = 1$, and using the classical convexity bound $L(s,\psi) \ll \eps^{-1} (\cond{\psi}(1+|t|))^{1+\eps-\sigma}$ for $-\eps \le \sigma \le 1+ \eps$, the horizontal and the vertical integrals are 
\eqs{ \ll W^\eps \prod_{\pi \mid ad} (1+\N(\pi)^{-\eps}) \bigl(W  T^{-1} + W^{\eps} T^{2- 2\eps}\bigr).}
Choosing $T =  W^{\frac{1-\eps}{3-2\eps}}$ shows that the double sum over $c_1, c_2$ equals
\eqs{R + O \Bigl(  \prod_{\pi \mid ad} (1+\N(\pi)^{-\eps})  W^{2/3 + 2 \eps} \Bigr),}
where the residue $R$ at $s=1$ is given by 
\bs{
	R &= \lim_{s\to 1} (W^s s^{-1} G(s) f^2 (s) )'\\
	&= A G(1) W \log W    -  A G(1) W + A G'(1)  W  + B G(1) W \\
	&= A G(1) W \log X -  W C. 
}
Here we wrote $G$ for $G_1$ and
\bs{f (s) &= (s-1) L(s,1) = (s-1) (1-3^{-s})\zeta_K (s), \\
	A &= \lim_{s\to 1} f^2 (s) = \frac{4\pi}{27\sqrt{3}} , \qquad B = \lim_{s\to 1}  (f^2 (s))' \\
	C&= A G(1) \log \N(de_1e_2)^2 +  A G(1)  - A G'(1)    - B G(1) .}
Therefore, we have shown so far that
\bs{\su{c \in \sF(X)\\ (c,a)=1} 1 &= 
\frac 1 {h_{(9)}} \su{d \equiv 1 \mod 3\\(d, a)=1\\ \N(d)^2 \le X} \mu_K(d)  \su{e_1 \equiv 1\mod 3\\ (e_1,ad)=1\\ \N(de_1)^2 \le X} \mu_K (e_1)  \su{e_2 \equiv 1\mod 3\\ (e_2, ad)=1\\ \N(de_1e_2)^2 \le X} \mu_K (e_2) \\
&\quad \cdot \Bigl( A G(1) W \log X -  W C + O \Bigl(  \prod_{\pi \mid ad} (1+\N(\pi)^{-\eps})  W^{2/3 + 2 \eps} \Bigr)  \Bigr).
}
Completing the sums over $d, e_1, e_2$ introduces an error $\ll \log \N(a) X^{1/2} \log^3 X$ since 
\eqn{\label{GG'}
	G'(1) \ll \log \N(ad), \qquad G(1) \le 1,}
and we obtain
\bs{\su{c \in \sF(X)\\ (c,a)=1} 1 
& = \frac {AX\log X} {h_{(9)}} \prod_{\pi \equiv 1 \mod 3} \Bigl(1 - \frac 3  {\N(\pi)^2} + \frac 2 {\N(\pi)^3} \Bigr) \pr{\pi \equiv 1 \mod 3 \\ \pi \mid a} \frac {\N(\pi)} {\N(\pi) + 2} \\
& + E(a) X +	O\Bigl( \prod_{\pi \mid a} (1+\N(\pi)^{-\eps}) X ^{2/3 + 2\eps} \Bigr), \\
}
where
\bs{
	E(a) & = -\frac 1 {h_{(9)}} \su{d \equiv 1 \mod 3\\(d, a)=1} \frac{\mu_K(d)}{\N(d)^2}  \su{e_1 \equiv 1\mod 3\\ (e_1,ad)=1} \frac{\mu_K (e_1)}{\N(e_1)^2} \su{e_2 \equiv 1\mod 3\\ (e_2, ad)=1}  \frac {\mu_K (e_2)} {\N(e_2)^2} C \\
& \ll \log \N(a), 
}
where inequality follows by \eqref{GG'}. Now recalling that $a = \cube$  it follows that
\eqn{\label{S1}
S_1 =  DX\log X + EX + O \bigl( X^{1+\eps} Y^{-1/6} + X ^{2/3 + 2\eps} \bigr),
}
where
\bsc{\label{AsymptoticConstant}
D &= \frac{4\pi}{81(\sqrt{3}-1)} \\
&\qquad \cdot \prod_{\pi \equiv 1 \mod 3} \Bigl(1 - \frac 3  {\N(\pi)^2} + \frac 2 {\N(\pi)^3} \Bigr) \Bigl(1 + \frac {\N(\pi)} {(\N(\pi) + 2)(\N(\xi)^{3/2}-1)} \Bigr), \\
E &= \frac 1 {1-1/\sqrt 3} \su{a \equiv 1 \mod 3\\ a = \cube} \frac {E(a)} {\N(a)^{1/2}}
.}
\subsection{Estimate of $S_2$}
Note that
\bsc{\label{S2integral}
&\su{a \equiv 1 \mod 3\\ a \neq \cube} \frac 1 {\N(a)^{1/2}} V \biggl( \frac {3^{r}\N(a)} Y \biggl) \su{c \in \sF(X)\\ (c,a)=1} \chi_{a} (c) \\
&= - \int_{1^-}^\infty  \biggl( \su{a \equiv 1 \mod 3\\ a \neq \cube\\ \N(a) \le z} \frac 1 {\N(a)^{1/2}} \su{c \in \sF(X)\\ (c,a)=1} \chi_{a} (c) \biggr)
\frac{3^r} Y V' \biggl( \frac {3^r z} Y \biggl) dz.
}
First consider the sums
\bs{
\su{c \in \sF\\ \cond c \sim y} \chi_{a} (c) &= \frac 1 {h_{(9)}} \sum_{\psi \mod 9} \su{d \equiv 1 \mod 3\\(d, a)=1\\ \N(d) \le y^{1/2}
} \mu_K(d)  \psi (d^3) \\
& \quad \cdot
\ps{c_1 \equiv 1 \mod 3\\  (c_1, ad)=1} (\chi_a \psi) (c_1)
\ps{c_2 \equiv 1 \mod 3\\ (c_2, ad)=1\\ \N(c_1c_2) \sim W} (\chi_a \psi) (c_2^2)
}
for $y=X2^{-k} \ge 1$ with $k \ge 1$, where $W = y/\N(d)^2$. It will be enough to estimate the sums
\eqn{\label{S2dyadic}
\su{a \equiv 1 \mod 3\\ (a,d) = 1 \\ a \neq \cube\\ \N(a) \sim z}   \frac  1  {\N(a)^{1/2}} \ps{c_1 \equiv 1 \mod 3\\ (c_1,ad) = 1\\  \N(c_1) \sim U} (\chi_a \psi) (c_1)
\ps{c_2 \equiv 1 \mod 3\\ \N(c_2) \sim V \\ (c_2,ad) = 1\\ \N(c_1c_2) \sim W} (\chi_a \psi) (c_2^2)
}
for a fixed $d$ with $\N(d) \in [1,\sqrt y]$ and $1 \le U = W2^{-i}, V=W2^{-j}$ with $i, j \ge 1$ satisfying $W/4 < UV < 2W$. Applying Perron's formula with $T = W$ for the sums over $c_1$ and $c_2$ and summing over $a$, we see that \eqref{S2dyadic} is bounded by 
\eqn{
	\label{S2AfterPerron}
	W \log W \sup_t	\su{a \equiv 1 \mod 3\\ (a,d) = 1 \\ a \neq \cube \\  \N(a) \sim z}   \frac 1 {\N(a)^{1/2}}  \big|\Sigma_1 (U) \Sigma_2 (V) \big|
	+ O\bigl( W^\eps z^{1/2} \bigr),
}
where 
\eqs{
	\Sigma_1 (U) = \ps{c_1 \equiv 1 \mod 3\\ (c_1,ad) = 1\\ \N(c_1) \sim U} \frac{(\chi_a \psi) (c_1)}{\N(c_1)^s}, \qquad 
	\Sigma_2 (V) = \ps{c_2 \equiv 1 \mod 3\\ (c_2,ad) = 1\\ \N(c_2) \sim V} \frac{(\chi_a \psi) (c_2^2)}{\N(c_2)^s},
}
and $s= 1 + 1/\log 2W + it$. Without loss of generality, we can assume that $U \le V$ so that $V \gg W^{1/2}$. Removing the square-free condition on $c_2$, we have that
\bs{\Sigma_2 (V) &= \su{e \equiv 1 \mod 3\\ (e,ad) = 1\\  \N(e) \le V^{1/2} } \frac{\mu_K(e) (\chi_a \psi) (e^4)}{\N(e)^{2s}}  \su{c_2 \equiv 1 \mod 3\\ (c_2,ad) = 1\\ \N(c_2e^2) \sim V} \frac{(\chi_a \psi) (c_2^2)}{\N(c_2)^s}.
}
For the inner sum, since $a \neq \cube$, we use the estimate
\eqs{\su{c_2 \equiv 1 \mod 3\\ (c_2,ad) = 1\\ \N(c_2e^2) \sim V} \frac{(\chi_a \psi) (c_2^2)}{\N(c_2)^s} \ll (V/\N(e)^2)^{-1/2+\eps} \N(a)^{1/4} \N(ad)^\eps,
}
which follows easily using Perron's formula and classical convexity bound together with partial integration, so that
\eqs{
\Sigma_2 (V) \ll V^{-1/2+\eps} \N(a)^{1/4} \N(ad)^\eps. 
}
Using this result in \eqref{S2AfterPerron} and applying cubic large sieve inequality for $\Sigma_1(U)$ yields that
\bs{
\eqref{S2AfterPerron} 
& \ll \N(d)^\eps \bigl(W^{ 1/2+\eps} z^{3/4+3\eps/2}  +  W^{3/4+\eps} z^{1/4+2\eps} +  W^{2/3+\eps}  z^{7/12+\eps}\bigr). 
}
Summing this over all $\ll \log W$ possible $U$'s and then over $\N(d) \le y^{1/2}$, and finally over $\ll \log X$ possible $y$'s, and inserting the resulting estimates in \eqref{S2integral} we conclude that 
\eqn{\label{S2bound}
S_2 \ll X^{ 1/2+\eps} Y^{3/4+3\eps/2}  +  X^{3/4+\eps} Y^{1/4+2\eps} + X^{2/3+\eps}  Y^{7/12+\eps}. 
}

\subsection{Estimate of $S_3$}
The contribution of $c$ with $\cond c \le Y$ to $S_3$ is $\ll Y\log Y$. Hence it is enough to consider 
\eqs{\sum_{a \equiv 1 \mod 3} \frac 1 {\N(a)^{1/2}} 
	\su{c \in \sF\\ Y < \cond c \le X}  \frac{\chi_c (a^2) W(\chi_c)}{\sqrt{\cond{c}}} V \biggl( \frac {Y\N(a)} {\cond{c}} \biggr).}
Using partial integration twice expresses this sum as
\bsc{\label{S3integral}
	& -  \int_{1^-}^\infty \frac Y X V' \biggl( \frac {Yz} X \biggr) \su{a \equiv 1 \mod 3\\ \N(a) \le z}   \frac 1 {\N(a)^{1/2}}  \su{c \in \sF\\ Y < \cond c \le X} \chi_c (a^2) \frac {W(\chi_c)}{\sqrt{\cond{c}}}  
	dz \\
	&- Y \int_{1^-}^\infty \int_Y^X \biggl( \frac 1 {y^2} V' \biggl( \frac {Yz} y \biggr) + \frac {Yz} {y^3} V'' \biggl( \frac {Yz} y \biggr) \biggr) \\
& \qquad \cdot \su{a \equiv 1 \mod 3\\ \N(a) \le z}   \frac 1 {\N(a)^{1/2}}  \su{c \in \sF\\ Y < \cond c \le y} \chi_c (a^2) \frac {W(\chi_c)}{\sqrt{\cond{c}}}  dy dz.
}
First note that for $c = c_1^2c_2 \in \sF$, 
\eqs{W(\chi_c) = \chi_{c_1^2} (c_2) \chi_{c_2} (c_1) \overline{W(\chi_{c_1})} W(\chi_{c_2}) = \overline{W(\chi_{c_1})} W(\chi_{c_2}),}
where the second equality follows by cubic reciprocity law. Note also that 
\eqs{W(\chi_{c_i}) = \chi_{c_i} (\sqrt{-3}) g(1,c_i),}
where $g(r,d)$ is the shifted Gauss sum defined in \eqref{ShiftedGaussSum}. Hence, 
\eqs{W(\chi_c) =  \chi_c (\sqrt{-3}) \overline{g(1,c_1)} g(1,c_2) = \overline{g(1,c_1)} g(1,c_2),}
since $\chi_c (\sqrt{-3}) = \chi_{c} (w(1-\omega)) = 1$ whenever $c \equiv 1 \mod 9$. Since $(a,c)=1$, it follows from \cite[Lemma 2.7]{DG} that $\chi_{c_1^2}(a^2)\overline{g(1, c_1)} = \overline{g(a,c_1)}$, and $\chi_{c_2}(a^2) g(1,c_2) = g(a,c_2)$. We can also remove the square-free conditions on $c_1, c_2$ using \cite[Lemma 2.8]{DG}. Thus, using ray class character modulo 9 to remove the condition $c_1^2c_2 \equiv 1 \mod 9$, we derive that 
\bs{\su{c \in \sF \cup \{1\} \\ \cond c \sim y} \chi_c (a^2) \frac{W(\chi_c)}{\sqrt{\cond{c}}} 
&= \frac 1 {h_{(9)}} \sum_{\psi \mod 9}
\su{c_1 \equiv 1 \mod 3\\ (c_1, a)=1} \psi (c_1^2) \frac{\overline{g(a, c_1)}}{\N(c_1)^{1/2}} \\
& \qquad \cdot 
\su{c_2 \equiv 1 \mod 3\\ (c_2,a c_1)=1 \\ \N(c_1c_2) \sim y/\N(d)^2} \psi (c_2) \frac{g(a,c_2)}{\N(c_2)^{1/2}}.
}
Finally, we remove the condition $(c_1, c_2)=1$ to get
\bs{
	&\su{c \in \sF \cup \{1\} \\ \cond c \sim y} \chi_c (a^2) \frac{W(\chi_c)}{\sqrt{\cond{c}}} 
	= \frac 1 {h_{(9)}} \sum_{\psi \mod 9} \su{d \equiv 1 \mod 3\\(d, a)=1\\ \N(d) \le \sqrt y} \mu_K(d) \psi (d^3)  \\
	&\qquad  \cdot  \su{c_1 \equiv 1 \mod 3\\  (c_1, ad)=1} \psi (c_1^2) \frac{\overline{g(ad, c_1)}}{\N(c_1)^{1/2}} 
	\su{c_2 \equiv 1 \mod 3\\ (c_2, ad)=1\\ \N(c_1c_2) \sim y/\N(d)^2} \psi(c_2) \frac{g(ad,c_2)}{\N(c_2)^{1/2}}
}
Here, we used \cite[Lemma 2.7]{DG} to first introduce the conditions $(d,c_1c_2)=1$, which can be done since otherwise the Gauss sums vanish, and then write $g(a,dc_i)$ as $g(ad, c_i)g(a,d)$, and finally used the fact that $|g(a,d)|^2 = \N(d)$ as $(d,a)=1$. 

We shall first estimate
\eqn{\label{S3dyadicaverage}
\su{a \equiv 1 \mod 3\\ (a,d) = 1 \\ \N(a) \sim z}   \frac 1 {\N(a)^{1/2}}  \su{c_1 \equiv 1 \mod 3\\  (c_1, ad)=1\\ \N(c_1) \sim U} \psi (c_1^2) \frac{\overline{g(ad, c_1)}}{\N(c_1)^{1/2}} 
\su{c_2 \equiv 1 \mod 3\\ (c_2, ad)=1\\ \N(c_2) \sim V \\ \N(c_1c_2) \sim W} \psi(c_2) \frac{g(ad,c_2)}{\N(c_2)^{1/2}}
}
for a fixed $d$ with $\N(d) \in [1,\sqrt y]$ and $1 \le U = W2^{-i}, V= W2^{-j}$ satisfying $W/4 < UV < 2W$, where $W= y/\N(d)^2$. Applying Perron's formula with $T = W$, we see that \eqref{S3dyadicaverage} is bounded by 
\eqn{
\label{LargeD}
W \log W \sup_t	\su{a \equiv 1 \mod 3\\ (a,d) = 1 \\ \N(a) \sim z}   \frac 1 {\N(a)^{1/2}}  \big|\Sigma_1 (U) \Sigma_2 (V) \big|
+ O\bigl( W^\eps z^{1/2} \bigr)
}
where 
\eqs{
\Sigma_1 (U) = \su{c_1 \equiv 1 \mod 3\\  (c_1, ad)=1\\ \N(c_1) \sim U} \psi (c_1^2) \frac{\overline{g(ad, c_1)}}{\N(c_1)^{1/2+s}}, \quad 
\Sigma_2 (V) = \su{c_2 \equiv 1 \mod 3\\ (c_2, ad)=1\\ \N(c_2) \sim V} \psi(c_2) \frac{g(ad,c_2)}{\N(c_2)^{1/2+s}},
}
and $s= 1 + 1/\log 2W + it$. 

By the last equation in the proof of \cite[Proposition 6.2]{DG},
\mult{\su{c \equiv 1 \mod 3\\  (c, r)=1\\ \N(c) \le z} \lambda(c) \frac{g(r, c)}{\N(c)^{1/2}} \\
\ll \N(r_1r_3^\ast)^\eps \bigl( z^{5/6}\N(r_1)^{-1/6} + z^{2/3+\eps} \N(r_1r_2^2)^{1/6} + z^{1/2+\eps} \N(r_1r_2^2)^{1/4} \bigr),	}
where $r=r_1r_2^2r_3^3$ with $r_i \equiv 1 \mod 3$, $r_1, r_2$ coprime and square-free and $r_3^\ast$ is the product of primes dividing $r_3$ but not $r_1r_2$. Also, the first error term appears only when $r_2 = 1$. Write $a = a_1a_2^2 a_3^3$ with $a_i \equiv 1\mod 3$, and $a_1, a_2$ square-free and coprime. Then, $ad = r_1r_2^2r_3^3$ with $r_1 = a_1d, r_2 = a_2$ and $r_3 = a_3$, since $(a,d)=1$ and $d$ is square-free. Using partial integration and the above result, we derive that 
\eqn{\label{PV}
	\Sigma_2 (V) \ll \N(a_1da_3)^\eps \Bigl(\frac 1 {V^{1/6} \N(a_1d)^{1/6}} +  \frac{\N(da_1a_2^2)^{1/6}}{V^{1/3-\eps}}  + \frac{\N(da_1a_2^2)^{1/4}}{V^{1/2-\eps}}  \Bigr) } 
where the first term appears only when $a_2 = 1$.  Using \eqref{PV} it then follows that the sum over $a$ in \eqref{LargeD} is bounded by 
\bs{
& \frac{(z\N(d))^\eps}{(z\N(d)V)^{1/6}}  \su{a_3 \equiv 1 \mod 3\\ \N(a_3) \le z^{1/3}} \N(a_3)^{-1} \biggl( \;  \ps{a_1\equiv 1 \mod 3\\ \N(a_1a_3^3) \sim z} \big| \Sigma_1 (U) \big|^2 \biggr)^{1/2}  
+  \frac{(z\N(d))^{1/6+\eps}  }{V^{1/3-\eps}} \\
& 
\cdot \sum_{\N(a_3)  \le z^{1/3}} \N(a_3)^{-2-2\eps} \ps{\N(a_2^2a_3^3) \le z} \frac 1 {\N(a_2)^{1+2\eps}} \biggl( \ps{a_1 \equiv 1 \mod 3\\ \N(a_1) \sim z/\N(a_2^2a_3^3)} |\Sigma_1(U)|^2\biggr)^{1/2} \\
& + \frac{(z\N(d))^{1/4+\eps} }{V^{1/2-\eps}} \sum_{\N(a_3)  \le z^{1/3}} \frac 1 {\N(a_3)^{9/4+2\eps}} \ps{\N(a_2^2a_3^3) \le z} \frac 1 {\N(a_2)^{1+2\eps}} \biggl( \ps{a_1 \equiv 1 \mod 3\\ \N(a_1) \sim z/\N(a_2^2a_3^3)} |\Sigma_1(U)|^2\biggr)^{1/2}.}
Applying cubic large-sieve inequality for each term and assuming that 
$U \le V$ shows that \eqref{LargeD} is bounded by 
\bs{
	& W^{11/12+ \eps} z^{-1/6+3\eps/2} \N(d)^{-1/6+\eps}
	+ W^{5/6+\eps} z^{1/3+3\eps/2} \N(d)^{-1/6+\eps}
	\\
	& +  W^{3/4+2\eps} \N(d)^{1/6+\eps} z^{1/2+3\eps/2}  
	+  W^{3/4+2\eps} \N(d)^{ 1/4+\eps} z^{1/4+3\eps/2} \\
	&
	+ W^{2/3+\eps} \N(d)^{1/6+\eps} z^{2/3+3\eps/2} 
	+  W^{2/3+2\eps} \N(d)^{ 1/4+\eps} z^{7/12+3\eps/2}
	\\
	& +  W^{1/2+2\eps} \N(d)^{ 1/4+\eps} z^{3/4+3\eps/2}    
	 +   W^{5/6+2\eps} \N(d)^{1/6+\eps} z^{1/6+3\eps/2}    + W^\eps z^{1/2}  .
}
Summing over the range $\N(d) \le B$ yields the bound
\bsc{\label{S3Smalld}
	& y^{11/12+ \eps} {z^{-1/6+3\eps/2} }
	+ y^{5/6+\eps} z^{1/3+3\eps/2} 
	+  y^{3/4+2\eps}  z^{1/2+3\eps/2}  
	\\
	&
	+ y^{2/3+\eps} z^{2/3+3\eps/2} 
	+  y^{1/2+2\eps} B^{ 1/4-\eps} z^{3/4+3\eps/2}   
	+ y^\eps B^{1-2\eps} z^{1/2}  .
}

Assume now that $B < \N(d) \le y^{1/2}$. Using Cauchy-Schwarz inequality gives
\eqs{
	\su{a \equiv 1 \mod 3\\ (a,d) = 1 \\ \N(a) \sim z}   \frac {\big|\Sigma_1 (U) \Sigma_2 (V) \big|} {\N(a)^{1/2}}   
	 \le z^{-1/2}\su{a_2 \equiv 1 \mod 3\\ \N(a_2) \le \sqrt z} 	\biggl( \; \ps{a_1 \equiv 1 \mod 3\\ \N(a_1a_2^2) \sim z}    \big|\Sigma_1 (U) \big|^2 \ps{a_1 \equiv 1 \mod 3\\ \N(a_1a_2^2) \sim z}   \big| \Sigma_2 (V) \big|^2 \biggr)^{1/2}.
}
By the cubic large sieve inequality,
\bs{
	\ps{a_1 \equiv 1 \mod 3\\ \N(a_1a_2^2) \sim z}    \big|\Sigma_1 (U) \big|^2 &\ll
	\Bigl(\frac z {U\N(a_2)^2}\Bigr)^{1+\eps}   +  \Bigl(\frac {Uz} {\N(a_2)^2}\Bigr)^{\eps}   + U^{-1/3+\eps} \Bigl(\frac z {\N(a_2)^2}\Bigr)^{2/3+\eps}, }
so that
\mult{
\su{a \equiv 1 \mod 3\\ (a,d) = 1 \\ \N(a) \sim z}   \frac {\big|\Sigma_1 (U) \Sigma_2 (V) \big|} {\N(a)^{1/2}}    \ll W^{\eps/2} \biggl(
	W^{-1/2} z^{1/2+\eps}   +  U^{-1/2} z^{\eps}   \\
+   U^{-1/2} z^{1/3+\eps} V^{-1/6}  +   V^{-1/2} z^{\eps}  +   1  +   z^{1/3+ \eps} U^{-1/6}  V^{-1/2} +   W^{-1/6} z^{1/6+\eps}  
	\biggr).
}
Inserting this in \eqref{LargeD} and
assuming $V \le U$ shows that \eqref{S3dyadicaverage} is bounded by
\eqs{
	W^{1/2+2\eps} z^{1/2+\eps}   +  W^{3/4+2\eps} z^{\eps}  +  W^{2/3+2\eps} z^{1/3+\eps} +  W^{1+\eps} z^{\eps}  +   z^{1/3+ \eps} W^{5/6+\eps} .
}
Summing over $B < \N(d) \le \sqrt y$ we get 	
\multn{\label{S3Larged}
y^{1/2+2\eps} z^{1/2+\eps}   +  y^{3/4+2\eps} z^{\eps} B^{-1/2} +  y^{2/3+2\eps} z^{1/3+\eps} B^{-1/3} \\
+  y^{1+\eps} z^{\eps} B^{-1} +    y^{5/6+\eps}  z^{1/3+ \eps} B^{-2/3}.
}

Combining \eqref{S3Smalld} and \eqref{S3Larged} and balancing terms using Lemma \ref{balancing} with $B \in [1,y^{1/2}]$ shows that for $y \in (Y,X]$, 
\eqs{ \su{a \equiv 1 \mod 3\\ \N(a) \le z}   \frac 1 {\N(a)^{1/2}}  \su{c \in \sF\\ Y < \cond c \le y} \chi_c (a^2) \frac {W(\chi_c)}{\sqrt{\cond{c}}}}
is bounded by
\mult{
y^{11/12+ \eps} z^{-1/6+3\eps/2} + y^{5/6+\eps} z^{1/3+3\eps/2} +  y^{3/4+2\eps}  z^{1/2+3\eps/2} 
\\
+ y^{2/3+\eps} z^{2/3+3\eps/2}  +  y^{1/2+2\eps} z^{3/4+3\eps/2}     +  y^{3/5+2\eps}  z^{3/5+2\eps} .
}
Using this result in \eqref{S3integral} we derive that
\eqn{\label{S3bound}
	S_3 
\ll X^\eps \bigl(X^{11/12} + X^{7/6} Y^{-1/3} 
+  X^{5/4}Y^{-1/2}    
 + X^{4/3} Y^{-2/3}  +  X^{6/5}Y^{-3/5} \bigr),
}
upon redefining $\eps$.
\subsection{Gluing the pieces}
Gathering \eqref{S1}, \eqref{S2bound} and \eqref{S3bound} we have shown that 
\bs{\sum_{c \in \sF(X)}& L(1/2, \chi_c)  
- DX\log X + EX \\
& \ll X^\eps \Bigl(X^{11/12} + X^{1/2} Y^{3/4}  +  X^{3/4} Y^{1/4} + \red{X^{2/3}  Y^{7/12}} \\
& \quad + X Y^{-1/6} + \red{X^{7/6} Y^{-1/3}} +  X^{5/4}Y^{-1/2} + X^{4/3} Y^{-2/3} +  X^{6/5}Y^{-3/5}\Bigr),  }
upon redefining $\eps$. 
Using Lemma \ref{balancing} with $Y \in [1,X]$ gives the claimed error term in theorem \eqref{thm:moments} and finishes the proof.

\begin{rem}
In the study conducted by Luo in \cite{Luo}, the first and second smoothed moments were examined within a thin subfamily of all cubic Dirichlet characters. In order to estimate the corresponding sum $S_2$, Luo utilized the cubic large sieve inequality for large values of $d$, along with a smoothed version of the Polya-Vinogradov inequality for small values of $d$. However, if one considers the first moment without using smoothing, as we have done in our work, it is sufficient to rely solely on the Polya-Vinogradov type inequality for the entire range of $1 \le \N(d) \le \sqrt y$.

Regarding the calculation of $S_3$, Luo employed a result concerning the shifted cubic sums, attributed to Patterson and referenced in \cite[equation 7]{Luo}. However, we encountered difficulty in locating this particular estimate in the cited paper. Furthermore, the mentioned estimate does not provide any indication of the dependence on the conductor of the character $\chi_a$, the shift parameter of the Gauss sums. In our previous work \cite{DG}, we addressed this issue. Upon reexamining Luo's case without incorporating the smooth sums, we discovered that it is still possible to obtain the first moment for the thin family, albeit with a slightly worse error term.
\end{rem}

\nocite{*}

\bibliographystyle{amsplain}

\end{document}